\documentclass[]{article}
\begin{document}
\newtheorem{proposition}{Proposition}[section]
\newtheorem{definition}{Definition}[section]
\newtheorem{lemma}{Lemma}[section]

\title{\bf Generalized Cohn's Theorem}
\author{Keqin Liu\\Department of Mathematics\\The University of British Columbia\\Vancouver, BC\\
Canada, V6T 1Z2}
\date{June 29, 2007}
\maketitle

\begin{abstract} We introduce the notion of a free associative $\mathcal{Z}_2$-algebra on the union of two disjoint sets and prove a generalization of Cohn's Theorem on Jordan algebras.\end{abstract}

It is well-known that the ordinary passage from an associative algebra $\mathcal{A}$ to a Lie algebra is established by the following binary operation:
\begin{equation}\label{eqa}
[x, y]=xy-yx \qquad\mbox{for $x$, $y\in \mathcal{A}$.}
\end{equation}
Replacing the sign $-$ by the sign $+$ and inserting a scalar factor $\displaystyle\frac12$ on the right hand side of (\ref{eqa}), the resulting binary operation:
\begin{equation}\label{eqb}
x\circ y=\displaystyle\frac12 (xy+yx) \qquad\mbox{for $x$, $y\in \mathcal{A}$}
\end{equation}
gives rise to the passage from an associative algebra to a Jordan algebra. 

\medskip
As a generalization of Lie algebras, (right) Leibniz algebras have been studied recently by a group of researchers. If $A$ is an associative $\mathcal{Z}_2$-algebra, then the binary operation:
\begin{equation}\label{eqc}
\langle x, y\rangle=xy_0-y_0x \qquad\mbox{for $x$, $y\in A$}
\end{equation}
gives a passage from an associative $\mathcal{Z}_2$-algebra to a (right) Leibniz algebra, where $y_0$ is the even component of $y$. Replacing the sign $-$ by the sign $+$ and inserting a scalar factor $\displaystyle\frac12$ on the right hand side of (\ref{eqc}), the resulting binary operation:
\begin{equation}\label{eqd}
x\bullet y=\displaystyle\frac12 (xy_0+y_0x) \qquad\mbox{for $x$, $y\in A$}
\end{equation}
introduces a class of non-associative algebras $(A, \, \bullet)$. Since the way of producing the binary operation (\ref{eqd}) from the binary operation (\ref{eqc}) is the same as the way of producing the binary operation (\ref{eqb}) from the binary operation (\ref{eqa}), we name the class of 
non-associative algebras $(A, \, \bullet)$  the generalized Jordan algebras in \cite{Liu}. Thus, the connection between Leibniz algebras and generalized Jordan algebra extends the connection between Lie algebras and Jordan algebra. Hence, if Leibniz algebras should be explored, then the generalized Jordan algebras  should not be ignored.

\medskip
The purpose of this paper is to prove a generalization of Cohn's Theorem on Jordan algebras in the context of the generalized Jordan algebras. In Section 1, we discuss the basic properties of a generalized Jordan algebra and establish the passage from an associative $\mathcal{Z}_2$-algebra to a generalized Jordan algebra. In Section 2, we state the universal property of the free associative $\mathcal{Z}_2$-algebra on the union of two disjoint sets. In Section 3, we introduce odd tetrads and extend Cohn's Theorem in the context of the generalized Jordan algebras.

\medskip
In this paper, all vector spaces are vector spaces over fields of characteristic not 2 and 3, and all associative algebras have an identity.

\medskip
\section{Generalized Jordan Algebras}

We begin this section with the definition of a generalized Jordan algebra, which was introduced in \cite{Liu}.

\begin{definition}\label{def1.1} A vector space $J$ is called a {\bf generalized Jordan algebra} if there exists a binary operation $\bullet: J\times J\to J$ such that the following two properties hold.
\begin{description}
\item[(i)] The binary operation $\bullet$ is {\bf right commutative}; that is
\begin{equation}\label{eq0}
 x\bullet (y\bullet z)=x\bullet (z\bullet y)\quad\mbox{for $x$, $y$, $z\in J$}.
\end{equation}
\item[(ii)] The binary operation $\bullet$ satisfies the {\bf Jordan identity}:
\begin{equation}\label{eq1}
(y\bullet x)\bullet (x\bullet x)=(y\bullet (x\bullet x))\bullet x
\end{equation}
and the {\bf Hu-Liu identity}:
\begin{equation}\label{eq2}
x\bullet (y\bullet (x\bullet x))-(x\bullet y)\bullet (x\bullet x)=
2(x\bullet x)\bullet (y\bullet x)-2((x\bullet x)\bullet y)\bullet x,
\end{equation}
where $x$, $y\in J$.
\end{description}
\end{definition}

A generalized Jordan algebra $J$ is also denoted by $(J, +, \bullet)$, where the binary operation $\bullet$ is called the {\bf bullet product}. If there exists an element $1$ of a generalized Jordan algebra $(J, +, \bullet)$ such that $x\bullet 1=x$ for $x\in J$, then $J$ is said to be  {\bf right unital} and the element $1$ is called a {\bf right unit}.

\begin{definition}\label{def1.2} An associative algebra $A$ is called an {\bf associative $\mathcal{Z}_2$-algebra} if $A=A_0\oplus A_1$ (as vector spaces) and
\begin{equation}\label{eq3}
A_0A_0\subseteq A_0, \quad A_0A_1+A_1A_0\subseteq A_1 \quad\mbox{and}\quad A_1A_1=0.
\end{equation}
\end{definition}

If $A=A_0\oplus A_1$ is an associative $\mathcal{Z}_2$-algebra, then $A_0$ and $A_1$ are called the {\bf even part} and {\bf odd part} of $A$ respectively. An element $x$ of $A=A_0\oplus A_1$ can be written uniquely as $x=x_0+x_1$, where $x_0\in A_0$ and $x_1\in A_1$ are called the {\bf even component} and {\bf odd component} of $x$, respectively. 

\medskip
The following proposition establishes the passage from an associative $\mathcal{Z}_2$-algebra to a right unital generalized Jordan algebra.

\begin{proposition}\label{pr1.1} If $A=A_0\oplus A_1$ is an associative $\mathcal{Z}_2$-algebra, then $A$ becomes a right unital generalized Jordan algebra under the following bullet product\begin{equation}\label{eq4}
x\bullet y=\displaystyle\frac12 (xy_0+y_0x),
\end{equation}
where $x$, $y\in A$ and $y_0$ is the even component of $y$.
\end{proposition}

\medskip
\noindent
{\bf Proof} First, the identity of the associative $\mathcal{Z}_2$-algebra $A$ is clearly a right unit by (\ref{eq4}).

\medskip
Next, for $x$, $y$, $z\in A$, we have
\begin{eqnarray*}
&&4x\bullet (y\bullet z)=2x\bullet (yz_0+z_0y)=x(yz_0+z_0y)_0+(yz_0+z_0y)_0x\\
&=&x(y_0z_0+z_0y_0)+(y_0z_0+z_0y_0)x=x(zy_0+y_0z)_0+(zy_0+y_0z)_0x\\
&=&2x(z\bullet y)_0+2(z\bullet y)_0 x=4x\bullet (z\bullet y),
\end{eqnarray*}
which proves that the bullet product $\bullet$ defined by (\ref{eq4}) is right commutative.

\medskip
Thirdly, we have
\begin{eqnarray*}
&&4(y\bullet x)\bullet (x\bullet x)=2(yx_0+x_0y)\bullet (x\bullet x)\\
&=&(yx_0+x_0y)(x\bullet x)_0+(x\bullet x)_0(yx_0+x_0y)\\
&=&(yx_0+x_0y)x^2_0+x^2_0(yx_0+x_0y)\\
&=&(yx^2_0)x_0+x_0(yx^2_0)+(x^2_0y)x_0+x_0(x_0^2y)\\
&=&(yx^2_0+x^2_0y)x_0+x_0(yx^2_0+x^2_0y)\\
&=&(y(x\bullet x)_0+(x\bullet x)_0y)x_0+x_0(y(x\bullet x)_0+(x\bullet x)_0y)\\
&=&2(y\bullet (x\bullet x))x_0+2x_0(y\bullet (x\bullet x))=(y\bullet (x\bullet x))\bullet x,
\end{eqnarray*}
which proves that the bullet product $\bullet$ defined by (\ref{eq4}) satisfies the Jordan identity.

\medskip
Finally, we have
\begin{eqnarray}
&&4x\bullet (y\bullet (x\bullet x))-4(x\bullet y)\bullet (x\bullet x)\nonumber\\
&=&2x\bullet (y(x\bullet x)_0+(x\bullet x)_0 y)
-2(x\bullet y)(x\bullet x)_0-2(x\bullet x)_0(x\bullet y)\nonumber\\
&=&2x\bullet (yx^2_0+x^2_0 y)-2(x\bullet y)x^2_0-2x^2_0(x\bullet y)\nonumber\\
&=&x (yx^2_0+x^2_0 y)_0+(yx^2_0+x^2_0 y)_0x-(xy_0+y_0x)x^2_0-x^2_0(xy_0+y_0x)\nonumber\\
&=&\underbrace{xy_0x^2_0}_{1}+xx^2_0 y_0+y_0x^2_0x+\underbrace{x^2_0 y_0x}_2
-\underbrace{xy_0x^2_0}_1-y_0xx^2_0-x^2_0xy_0-\underbrace{x^2_0y_0x}_2\nonumber\\
\label{eq5}&=&xx^2_0 y_0+y_0x^2_0x-y_0xx^2_0-x^2_0xy_0
\end{eqnarray}
and
\begin{eqnarray}
&&8(x\bullet x)\bullet (y\bullet x)-8((x\bullet x)\bullet y)\bullet x\nonumber\\
&=&2(xx_0+x_0x)\bullet (yx_0+x_0y)-4((xx_0+x_0x)\bullet y)\bullet x\nonumber\\
&=&(xx_0+x_0x) (yx_0+x_0y)_0+(yx_0+x_0y)_0(xx_0+x_0x) +\nonumber\\
&&\quad -2((xx_0+x_0x) y_0+y_0(xx_0+x_0x))\bullet x\nonumber\\
&=&(xx_0+x_0x) (y_0x_0+x_0y_0)+(y_0x_0+x_0y_0)(xx_0+x_0x) +\nonumber\\
&&\quad -(xx_0y_0+x_0xy_0+y_0xx_0+y_0x_0x)x_0+\nonumber\\
&&\quad -x_0(xx_0y_0+x_0xy_0+y_0xx_0+y_0x_0x)\nonumber\\
&=&\underbrace{xx_0y_0x_0}_1+xx_0x_0y_0+\underbrace{x_0xy_0x_0}_2+\underbrace{x_0xx_0y_0}_3
+\underbrace{y_0x_0xx_0}_4+y_0x_0x_0x+\nonumber\\
&&\quad +\underbrace{x_0y_0xx_0}_5+\underbrace{x_0y_0x_0x}_6-\underbrace{xx_0y_0x_0}_1-\underbrace{x_0xy_0x_0}_2-y_0xx_0x_0+\nonumber\\
&&\quad -\underbrace{y_0x_0xx_0}_4-\underbrace{x_0xx_0y_0}_3-x_0x_0xy_0-\underbrace{x_0y_0xx_0}_5-\underbrace{x_0y_0x_0x}_6\nonumber\\
\label{eq6}&=&xx^2_0 y_0+y_0x^2_0x-y_0xx^2_0-x^2_0xy_0
\end{eqnarray}

It follows from (\ref{eq5}) and (\ref{eq6}) that the  the bullet product $\bullet$ defined by (\ref{eq4}) satisfies the Hu-Liu bullet identity.

\medskip
This completes the proof of Proposition~\ref{pr1.1}.

\hfill\raisebox{1mm}{\framebox[2mm]{}}

\bigskip
{\bf Remark} Except associative $\mathcal{Z}_2$-algebras, there exists another class of associative algebras which can be also used to generalize the well-known passage from an associative algebra to a Jordan algebra. The class of associative algebras are called invariant algebras in \cite{Liu3}. In fact, if we replace the sign $-$ appearing among the four Hu-Liu angle brackets introduced in Section 1.3 of \cite{Liu3} by the sign $+$, then resulting new binary operations give a few passages from an invariant algebra to a generalized Jordan algebra.

\medskip
For convenience, we will use $A^{\mathring{+}}$ to denote the right unital generalized Jordan algebra obtained from an associative $\mathcal{Z}_2$-algebra $A=A_0\oplus A_1$ by using the bullet product (\ref{eq4}). Note that $A_0$ is a ordinary Jordan algebra with respect to the bullet product (\ref{eq4}).

\begin{definition}\label{def1.3} Let $I$ be a subspace of a  generalized Jordan algebra 
$(J, +, \bullet)$.
\begin{description}
\item[(i)] I is called an {\bf ideal} of $J$ if $I\bullet J+J\bullet I\subseteq I$.
\item[(ii)] I is called a {\bf (generalized Jordan) subalgebra} of $J$ if 
$I\bullet  I\subseteq I$.
\end{description}
\end{definition}

The {\bf annihilator} $J^{ann}$ of a generalized Jordan algebra $(J, +, \bullet)$ over a vector space $\mathbf{k}$ is defined by
\begin{equation}\label{eq01}
J^{ann}:=\sum _{x, y\in J}\mathbf{k} (x\bullet y-y\bullet x).
\end{equation}

For $x$, $y$, $z\in J$, we have
$$
z\bullet (x\bullet y-y\bullet x)=z\bullet (x\bullet y)-z\bullet (y\bullet x)=0
$$
and
$$
(x\bullet y-y\bullet x)\bullet z=\Big((x\bullet y)\bullet z-z\bullet (x\bullet y)\Big)+\Big(z\bullet (y\bullet x)-(y\bullet x)\bullet z\Big)
$$
by the right commutative law. Hence, the annihilator $J^{ann}$ is an ideal of $J$.

\begin{definition}\label{def1.4} A generalized Jordan algebra  $J$ is called a {\bf simple generalized Jordan algebra} if $J\ne J^{ann}$, $J^{ann}\ne 0$ and $J$ has no ideals which are not equal to $\{0\}$, $J^{ann}$ and $J$.
\end{definition}

If $(J, +, \bullet)$ is a right unital generalized Jordan algebra, then
$$
J^{ann}=\{\, a \,|\, 1\bullet a=0\,\}
$$
and
$$
\{\, 1+a \,|\, a\in J^{ann}\,\}=\mbox{the set of all right units of $J$},
$$
where $1$ is a right unit of $J$.

\medskip
Let $(J, +, \bullet)$ be a generalized Jordan algebra. The {\bf long associator} $[x, y, z]_\ell$ is defined by
\begin{equation}\label{eq02}
[x, y, z]_\ell :=x\bullet (y\bullet z)-(x\bullet y)\bullet z-2z\bullet (y\bullet x)
+2(z\bullet y)\bullet x,
\end{equation}
where $x$, $y$, $z\in J$. By the Hu-Liu identity, we have
\begin{equation}\label{eq03}
[x, y, x\bullet x]_\ell =0 \quad\mbox{for $x$, $y\in J$.}
\end{equation}

After linearizing (\ref{eq03}), we get
\begin{equation}\label{eq04}
[x, y, x\bullet z]_\ell +[x, y, z\bullet x]_\ell+[z, y, x\bullet x]_\ell=0 
\end{equation}
and
\begin{equation}\label{eq05}
[x, y, w\bullet z+z\bullet w]_\ell + [w, y, z\bullet x+x\bullet z]_\ell 
+[z, y, x\bullet w+w\bullet x]_\ell=0 ,
\end{equation}
where $x$, $y$, $z$, $w\in J$. 

For $a\in J$, let $L_a: J\to J$ and $R_a: J\to J$ be the {\bf left multiplication} and the {\bf right multiplication}; that is
$$
L_a(x): =a\bullet x, \quad R_a(x): =x\bullet a \quad\mbox{for all $x\in J$.}
$$
For convenience, we also define $S_a$ by
$$
S_a:=L_a +R_a \quad\mbox{for all $x\in J$.}
$$

It follows from (\ref{eq0}), (\ref{eq1}) and (\ref{eq2}) that left multiplications and right multiplications have the following properties:
\begin{equation}\label{eq06}
L_xL_y=L_xR_y, \quad R_{x\bullet y}=R_{y\bullet x} \quad\mbox{for all $x$, $y\in J$,}
\end{equation}
\begin{equation}\label{eq07}
R_xR_{x\bullet x}=R_{x\bullet x}R_x \quad\mbox{for all $x\in J$}
\end{equation}
and
\begin{equation}\label{eq08}
L_xR_{x\bullet x}-R_{x\bullet x}L_x=2L_{x\bullet x}R_x-2R_xL_{x\bullet x} \quad\mbox{for all $x\in J$.}
\end{equation}
The long associator $[x, y, z]_\ell$ can also be written as
\begin{equation}\label{eq09}
[x, y, z]_\ell \left\{\begin{array}{l}
=(R_{y\bullet z}-R_zR_y-2L_zL_y+2L_{z\bullet y})(x),\\
=(L_xR_z-R_zL_x-2L_zR_x+2R_xL_z)(y),\\
=(L_xL_y-L_{x\bullet y}-2R_{y\bullet x}+2R_xR_y)(z).\end{array}\right.
\end{equation}

By (\ref{eq09}), (\ref{eq05}) is equivalent to
\begin{eqnarray}\label{eq010}
&&(L_xL_y-L_{x\bullet y}-2R_{y\bullet x}+2R_xR_y)S_z+\nonumber\\
&&\quad +2(R_{y\bullet (z\bullet x)}-R_{z\bullet x}R_y-L_{z\bullet x+x\bullet z}L_y
+L_{(z\bullet x+x\bullet z)\bullet y})+\nonumber\\
&&\quad +(L_zL_y-L_{z\bullet y}-2R_{y\bullet z}+2R_zR_y)S_x=0.
\end{eqnarray}

Letting $x=z$ in (\ref{eq010}), we get
\begin{eqnarray}\label{eq011}
&&(L_xL_y-L_{x\bullet y}-2R_{y\bullet x}+2R_xR_y)S_x+\nonumber\\
&&\quad 
+R_{(x\bullet x)\bullet y}-R_{x\bullet x}R_y-2L_{x\bullet x}L_y+2L_{(x\bullet x)\bullet y}=0.
\end{eqnarray}

\medskip
Recall from Section 9 of Chapter II in \cite{N} that a vector space $V$ over a Jordan algebra $(\mathcal{J}, +, \odot)$ is called a {\bf Jordan bimodule} if there is a bilinear map $(v, a)\mapsto va$ form $V\times\mathcal{J} \to V$ satisfying
\begin{equation}\label{eq012}
(v (a \odot a))a=(va)(a \odot a)
\end{equation}
and
\begin{equation}\label{eq013}
2((va)b)a+v((a \odot a) \odot b)=2(va)(a \odot b)+(vb)(a \odot a)
\end{equation}
for $v\in V$ and $a$, $b\in\mathcal{J}$.

\medskip
The next proposition shows that a generalized Jordan algebra structure is obtained by combining of a Jordan algebra and a Jordan bimodule.

\begin{proposition}\label{pr1.2} If $(J, +, \bullet)$ is a generalized Jordan algebra over a field $\mathbf{k}$, then the annihilator $J^{ann}$ becomes a bimodule over the Jordan algebra
$\displaystyle\frac{J}{J^{ann}}$ under the following bimodule action:
\begin{equation}\label{eq014}
u\bar{x}:= u\bullet x \quad\mbox{for $u\in J^{ann}$, $x\in J$ and $\bar{x}:=x+J^{ann}\in
\displaystyle\frac{J}{J^{ann}}$.}
\end{equation}
\end{proposition}

\medskip
\noindent
{\bf Proof} Note that $\left(\displaystyle\frac{J}{J^{ann}}, +, \odot\right)$ is a Jordan algebra, where the product $\odot$ is defined by 
$$\bar{x}\odot \bar{y}:=\overline{x\bullet y}\quad\mbox{for $x$, $y\in J$.}$$

\medskip
By the right commutative property, the action (\ref{eq014}) is well-defined. According to (\ref{eq012}) and (\ref{eq013}), we need to prove
\begin{equation}\label{eq015}
(u (\bar{x}\odot \bar{x}))\bar{x}=(u\bar{x})(\bar{x}\odot \bar{x})
\end{equation}
and
\begin{equation}\label{eq016}
2((u\bar{x})\bar{y})\bar{x}+u((\bar{x}\odot \bar{x})\odot \bar{y})
=2(u\bar{x})(\bar{x}\odot \bar{y})+(u\bar{y})(\bar{x}\odot \bar{x})
\end{equation}
for $u\in J^{ann}$ and $x$, $y\in J$.

\medskip
Since 
$$u (\bar{x}\odot \bar{x})=u\,\,\overline{x\bullet x}=u\bullet(x\bullet x)
=(R_{x\bullet x}|J^{ann})(u),$$
(\ref{eq015}) follows from (\ref{eq07}).

\medskip
Using $L_{x}|J^{ann}=0$ and $S_{x}|J^{ann}=R_{x}|J^{ann}$, we get from 
(\ref{eq011}) that
$$
-2R_{y\bullet x}R_x(u)+2R_xR_yR_x(u)=R_{(x\bullet x)\bullet y}(u)-R_{x\bullet x}R_y(u)
$$
or
$$
-2(u\bar{x})(\bar{x}\odot \bar{y})+2((u\bar{x})\bar{y})\bar{x}
=u((\bar{x}\odot \bar{x})\odot \bar{y})-(u\bar{y})(\bar{x}\odot \bar{x}),
$$
which is (\ref{eq016}).

\hfill\raisebox{1mm}{\framebox[2mm]{}}

We now give a way of constructing a generalized Jordan algebra structure from a Jordan algebra and a Jordan bimodule over the Jordan algebra. 

\begin{proposition}\label{pr1.3} If $(\mathcal{J}, +, \odot)$ is a Jordan algebra and $V$ a Jordan bimodule over the Jordan algebra $(\mathcal{J}, +, \odot)$, then the vector space direct sum 
$$J:=\mathcal{J}\oplus V=\{\, (a, v) \,|\, \mbox{$a\in\mathcal{J}$ and $v\in V$} \,\}$$
becomes a generalized Jordan algebra under the following bullet product:
\begin{equation}\label{eq017}
(a, v)\bullet (b, u):=(a\odot b, vb)\quad\mbox{$a$, $b\in\mathcal{J}$ and $v$, $u\in V$.}
\end{equation}
\end{proposition}

\medskip
\noindent
{\bf Proof} Let $(a, v)$, $(b, u)$, $(c, w)\in \mathcal{J}\oplus V$, where $a$, $b$, $c\in \mathcal{J}$ and $v$, $u$, $w\in V$.

\medskip
First, we have
\begin{eqnarray*}
&&(a, v)\odot \Big((b, u)\odot (c, w)\Big)=(a, v)\odot (b\odot c, uc)
=\Big(a\odot (b\odot c), v(b\odot c)\Big)\\
&=&\Big(a\odot (c\odot b), v(c\odot b)\Big)=(a, v)\odot (c\odot b, wb)
=(a, v)\odot \Big((c, w)\odot (b, u)\Big),
\end{eqnarray*}
which proves that the bullet product $\bullet$ defined by (\ref{eq017}) is right commutative.

\medskip
Next, we have
\begin{eqnarray}\label{eq018}
&&\Big((b, u)\bullet (a, v)\Big)\bullet \Big((a, v)\bullet (a, v)\Big)=
(b\odot a, ua)\bullet (a\odot a, va)\nonumber\\
&=&\Big((b\odot a)\odot (a\odot a), (ua)(a\odot a)\Big)
\end{eqnarray}
and
\begin{eqnarray}\label{eq019}
&&\Big((b, u)\bullet \big((a, v)\bullet (a, v)\big)\Big)\bullet (a, v)=
\Big((b, u)\bullet (a\odot a, va)\Big)\bullet (a, v)\nonumber\\
&=&\Big(b\odot (a\odot a), u(a\odot a)\Big)\bullet (a, v)\nonumber\\
&=&\Big(\big(b\odot (a\odot a)\big)\odot a, \big(u(a\odot a)\big)a\Big).
\end{eqnarray}
By (\ref{eq018}) and (\ref{eq019}), the bullet product $\bullet$ defined by (\ref{eq017}) satisfies the Jordan identity.

\medskip
Finally, we have
\begin{eqnarray}\label{eq020}
&&(a, v)\bullet \Big((b, u)\bullet \big((a, v)\bullet (a, v)\big)\Big)-
\Big((a, v)\bullet (b, u)\Big)\bullet \Big((a, v)\bullet (a, v)\Big)\nonumber\\
&=&(a, v)\bullet \Big((b, u)\bullet (a\odot a, va)\Big)-
\Big((a\odot b, vb)\bullet (a\odot a, va)\Big)\nonumber\\
&=&(a, v)\bullet \Big(b\odot (a\odot a), u(a\odot a)\Big)-
\Big((a\odot b)\odot (a\odot a), (vb)(a\odot a)\Big)\nonumber\\
&=&\Big(a\odot \big(b\odot (a\odot a)\big), v\big(b\odot (a\odot a)\big)\Big)-
\Big((a\odot b)\odot (a\odot a), (vb)(a\odot a)\Big)\nonumber\\
&=&\Big(0, v\big((a\odot a)\odot b\big)-(vb)(a\odot a)\Big)
\end{eqnarray}
and
\begin{eqnarray}\label{eq021}
&&2\Big(\big((a, v)\bullet (a, v)\big)\bullet \big((b, u)\bullet (a, v)\big)\Big)-
2\Big(\big((a, v)\bullet (a, v)\big)\bullet (b, u)\Big)\bullet (a, v)\nonumber\\
&=&2\Big((a\odot a, va)\bullet (b\odot a, ua)\Big)-
2\Big((a\odot a, va)\bullet (b, u)\Big)\bullet (a, v)\nonumber\\
&=&2\Big((a\odot a)\odot (b\odot a), (va)(b\odot a)\Big)-
2\Big((a\odot a)\odot b, (va)b\Big)\bullet (a, v)\nonumber\\
&=&2\Big((a\odot a)\odot (b\odot a), (va)(b\odot a)\Big)-
2\Big(\big((a\odot a)\odot b\big)\odot a, \big((va)b\big)a\Big)\nonumber\\
&=&\Big(0, 2(va)(a\odot b)\Big)-2\big((va)b\big)a\Big).
\end{eqnarray}
By (\ref{eq020}) and (\ref{eq021}), the bullet product $\bullet$ defined by (\ref{eq017}) satisfies the Hu-Liu identity.

\medskip
This completes the proof of Proposition~\ref{pr1.3}.

\hfill\raisebox{1mm}{\framebox[2mm]{}}

\medskip
\section{Free Associative $\mathcal{Z}_2$-Algebra}

Let $\check{X}$ and $\check{\Theta}$ be two disjoint sets. Let 
$\mathcal{F}\mathcal{A}[\check{X}\cup\check{\Theta}]$ denote the free unital $\mathbf{k}$-associative algebra on the set $\check{X}\cup\check{\Theta}$; that is, $\mathcal{F}\mathcal{A}[\check{X}\cup\check{\Theta}]$ is a vector space over $\mathbf{k}$ with a basis consisting of all {\bf monomials} $\check{u}_1\cdots \check{u}_n$ for all $n\ge 0$ ( the {\bf empty product} for $n=0$ serving as unit $1$) and all 
$\check{u}_i\in \check{X}\cup\check{\Theta}$, with the associative product determined by linearity and juxtaposition:
$$
(\check{u}_1\cdots \check{u}_n)(\check{u}_{n+1}\cdots \check{u}_{n+m})=
\check{u}_1\cdots \check{u}_n\check{u}_{n+1}\cdots \check{u}_{n+m} .
$$
The $\Theta$-degree $deg_{_\Theta}(\check{u}_1\cdots \check{u}_n)$ of a monomial 
$\check{u}_1\cdots \check{u}_n$ is defined by
$$
deg_{_\Theta}(\check{u}_1\cdots \check{u}_n):=\left\{\begin{array}{ll}
0&\mbox{if $n=0$},\\|\{\check{u}_i\,|\, \mbox{$\check{u}_i\in\Theta$ and $1\le i\le n$}\}|
&\mbox{if $n\ge 1$}.\end{array}\right.
$$
Let
$$
I_{_\Theta}:=\displaystyle\sum_{\begin{array}{l}\mbox{$\check{u}$ is a monomial}\\
\mbox{and $deg_{_\Theta}(\check{u})\ge 2$}\end{array}}\mathbf{k}\check{u} .
$$
Then $I_{_\Theta}$ is an ideal of the free associative algebra $\mathcal{F}\mathcal{A}[\check{X}\cup\check{\Theta}]$. We use 
$\mathcal{F}\mathcal{A}_2[X\cup\Theta]$ to denote the quotient associative algebra
$\displaystyle\frac{\mathcal{F}\mathcal{A}[\check{X}\cup\check{\Theta}]}{I_{_\Theta}}$, where
$$
X:=\{\, x \,|\, \mbox{$x=\check{x}+I_{_\Theta}$ and $\check{x}\in\check{X}$}\,\}
$$
and
$$
\Theta:=\{\, \theta \,|\, \mbox{$\theta =\check{\theta }+I_{_\Theta}$ and $\check{\theta }\in\check{\Theta }$}\,\}.
$$
Then $\mathcal{F}\mathcal{A}_2[X\cup\Theta]$ is an associative $\mathcal{Z}_2$-algebra whose even part $\mathcal{F}\mathcal{A}_2[X\cup\Theta]_0$ and odd part $\mathcal{F}\mathcal{A}_2[X\cup\Theta]_1$ are given by
$$
\mathcal{F}\mathcal{A}_2[X\cup\Theta]_0:=\displaystyle\bigoplus_{\begin{array}{c}n\ge 0\\
x_{i_1}, \cdots , x_{i_n}\in X\end{array}}\mathbf{k}x_{i_1}\cdots x_{i_n}
$$
and
$$
\mathcal{F}\mathcal{A}_2[X\cup\Theta]_1:=\displaystyle\bigoplus_{\begin{array}{c}n, m\ge 0\\
x_{i_1}, \cdots , x_{i_n}\in X\\ x_{j_1}, \cdots , x_{j_m}\in X\\
\theta_{x_{i_1}, \cdots , x_{i_n}\atop x_{j_1}, \cdots , x_{j_m}}\in\Theta \end{array}}
\mathbf{k}x_{i_1}\cdots x_{i_n}\theta_{x_{i_1}, \cdots , x_{i_n}\atop x_{j_1}, \cdots , x_{j_m}}
x_{j_1}\cdots  x_{j_m} .
$$
$\mathcal{F}\mathcal{A}_2[X\cup\Theta]$ is called the {\bf free associative 
$\mathcal{Z}_2$-algebra} on the set $X\cup\Theta$. The next proposition gives the universal property of $\mathcal{F}\mathcal{A}_2[X\cup\Theta]$.

\begin{proposition}\label{pr2.1} If $\phi$ is a map from the set $X\cup\Theta$ to an associative $\mathcal{Z}_2$-algebra $A=A_0\oplus A_1$ such that $\phi (\Theta )\subseteq A_1$, then $\phi$ can be extended uniquely to an associative algebra homomorphism from $\mathcal{F}\mathcal{A}_2[X\cup\Theta]$ to $A$; that is, there exists a unique associative algebra homomorphism $\tilde{\phi}: \mathcal{F}\mathcal{A}_2[X\cup\Theta] \to A$ such that
$\tilde{\phi}|(X\cup\Theta)=\phi$.
\end{proposition}

\medskip
\noindent
{\bf Proof} $\tilde{\phi}$ is clearly unique. We need only to prove the existence of $\tilde{\phi}$. By the universal property of the free unital associative algebra $\mathcal{F}\mathcal{A}[\check{X}\cup\check{\Theta}]$, there exists an associative algebra homomorphism $\check{\phi}: \mathcal{F}\mathcal{A}[\check{X}\cup\check{\Theta}]\to A$ such that
\begin{equation}\label{eq7}
\check{\phi}(\check{y})=\phi (y),
\end{equation}
where $y=\check{y}+\check{\Theta}$ and $\check{y}\in \check{X}\cup\check{\Theta}$. If $\check{u}$ is a monomial with $deg_{_\Theta}(\check{u})\ge 2$, then there exist 
$\check{y}_{i_1}, \cdots , \check{y}_{i_t}\in \check{X}\cup\check{\Theta}$ and 
$\check{\theta} _1$, $\check{\theta} _2\in\\check{Theta}$ such that
$$
\check{u}=\check{y}_{i_1} \cdots \check{y}_{i_n}\,\check{\theta} _1\,\check{y}_{i_{n+1}} \cdots \check{y}_{i_{n+m}}\,\check{\theta} _2\,\check{y}_{i_{n+m+1}} \cdots \check{y}_{i_t}.
$$
Since $\phi (\theta _i)\in A_1$ for $i=1$ and $2$, it follows from (\ref{eq7}) that
\begin{eqnarray*}
\check{\phi}(\check{u})&=&
\check{\phi}(\check{y}_{i_1}) \cdots \check{\phi}(\check{y}_{i_n})
\,\check{\phi}(\check{\theta} _1)\,\check{\phi}(\check{y}_{i_{n+1}}) \cdots \check{\phi}(\check{y}_{i_{n+m}})
\,\check{\phi}(\check{\theta} _2)\,\check{\phi}(\check{y}_{i_{n+m+1}}) \cdots \check{\phi}(\check{y}_{i_t})\\
&=&y_{i_1} \cdots y_{i_n}\,\theta _1\,y_{i_{n+1}} \cdots y_{i_{n+m}}
\,\theta _2\,y_{i_{n+m+1}} \cdots y_{i_t}=0,
\end{eqnarray*}
which proves that $I_{_\Theta}\subseteq Ker \check{\phi}$. Thus, $\check{\phi}$ induces an associative algebra homomorphism
$$
\tilde{\phi}: \mathcal{F}\mathcal{A}_2[X\cup\Theta]=
\displaystyle\frac{\mathcal{F}\mathcal{A}[\check{X}\cup\check{\Theta}]}{I_{_\Theta}}\to A
$$
such that $\tilde{\phi}(y)=\check{\phi}(\check{y})=\phi (y)$ for $y\in X\cup\Theta$. This proves Proposition~\ref{pr2.1}.

\hfill\raisebox{1mm}{\framebox[2mm]{}}

\medskip
\section{Generalized Cohn's Theorem}

Let $(\mathcal{F}\mathcal{A}_2[X\cup\Theta]^{\mathring{+}}, +, \bullet)$ be the right unital generalized Jordan algebra obtained from the free associative $\mathcal{Z}_2$-algebra $\mathcal{F}\mathcal{A}_2[X\cup\Theta]$, where the bullet product $\bullet$ is defined by (\ref{eq4}). The free associative $\mathcal{Z}_2$-algebra $\mathcal{F}\mathcal{A}_2[X\cup\Theta]$ has a unique {\bf reversal involution} $\ast$ such that
$$
1^\ast :=1, \quad (y_1y_2\cdots y_n)^\ast :=y_n^\ast\cdots y_2^\ast y_1^\ast \quad\mbox{for $y_1, \cdots , y_n\in X\cup\Theta$. }
$$
It is clear that the reversal involution $\ast$ preserves the $\mathcal{Z}_2$-grading of 
$\mathcal{F}\mathcal{A}_2[X\cup\Theta]$. An element $a$ of $\mathcal{F}\mathcal{A}_2[X\cup\Theta]$ is said to be {\bf reversible} if $a^\ast =a$. If $a$ and $b$ are two reversible elements of $\mathcal{F}\mathcal{A}_2[X\cup\Theta]$, then
\begin{eqnarray*}
&&(a\bullet b)^\ast =\left(\frac12(ab_0+b_0a)\right)^\ast 
=\frac12\Big((b_0)^\ast a^\ast +a^\ast (b_0)^\ast \Big)\\
&=&\frac12\Big((b^\ast)_0 a^\ast +a^\ast (b^\ast)_0 \Big)=\frac12(b_0a+ab_0)=a\bullet b.
\end{eqnarray*}
This proves that the reversible elements of $\mathcal{F}\mathcal{A}_2[X\cup\Theta]$ form a right unital generalized Jordan subalgebra $\mathcal{H}(\mathcal{F}\mathcal{A}_2[X\cup\Theta], \ast)$ of $\mathcal{F}\mathcal{A}_2[X\cup\Theta]^{\mathring{+}}$. 

\medskip
In the remaining of this section, we fix an ordering of $X\cup\Theta$ with the following property:
$$
\theta <x \quad\mbox{for $\theta\in \Theta$ and $x\in X$.}
$$
For $y_1, \cdots , y_n\in X\cup\Theta$, let
$$
\{y_1, y_2, \cdots , y_n\}:=\frac12\,(y_1 y_2 \cdots  y_n +y_n\cdots y_2 y_1).
$$
If $x_1$, $x_2$, $x_3$, $x_4$ are distinct elements of $X$ and $\theta\in \Theta$, then 
$\{x_1, x_2, x_3, x_4\}$ with $x_1< x_2< x_3< x_4$ is called an {\bf even tetrad}, and 
$\{\theta , x_1, x_2, x_3\}$ with $x_1< x_2< x_3$ is called an {\bf odd tetrad}. 

\medskip
The following proposition gives a generalization of Cohn's Theorem.

\begin{proposition}\label{pr3.1} The right unital generalized Jordan algebra  $\mathcal{H}(\mathcal{F}\mathcal{A}_2[X\cup\Theta], \ast)$ of reversible elements of the free associative $\mathcal{Z}_2$-algebra $\mathcal{F}\mathcal{A}_2[X\cup\Theta]$ coincides with the generalized Jordan subalgebra $\mathcal{H}'$ of $\mathcal{F}\mathcal{A}_2[X\cup\Theta]^{\mathring{+}}$ generated by $1$, $X\cup\Theta$ and all even tetrads and odd tetrads.
\end{proposition}

\medskip
\noindent
{\bf Proof} It is clear that $\mathcal{H}'\subseteq \mathcal{H}(\mathcal{F}\mathcal{A}_2[X\cup\Theta], \ast)$. Hence, we need only to prove
\begin{equation}\label{eq8}
\mathcal{H}(\mathcal{F}\mathcal{A}_2[X\cup\Theta], \ast)\subseteq \mathcal{H}' .
\end{equation}

\medskip
It is easy to check that $\mathcal{H}(\mathcal{F}\mathcal{A}_2[X\cup\Theta], \ast)$ is spanned by the set
$$
\left\{\, \{x_1, x_2, \cdots , x_n\}, \, 
\{x_1, x_2, \cdots , x_n, \theta , x_{n+1}, \cdots , x_m\}\,\left|\, 
\begin{array}{c} n,m\ge 0,\, \theta\in\Theta \\x_1, \cdots , x_m\in X \end{array}\right.\right\},
$$
where
$$
\{x_1, x_2, \cdots , x_n, \theta , x_{n+1}, \cdots , x_m\}:=\left\{\begin{array}{ll}
\{x_1, x_2, \cdots , x_n, \theta \}&\mbox{if $n=m$},\\
\{\theta , x_{1}, x_{2},\cdots , x_m\}&\mbox{if $n=0$}.\end{array}\right.
$$
Hence, in order to prove (\ref{eq8}), it is enough to prove
\begin{equation}\label{eq9}
\{x_1, x_2, \cdots , x_n\}\,\equiv \,0\,(\,\mbox{mod}\, \mathcal{H}') \quad\mbox{for $n\ge 0$}
\end{equation}
and
\begin{equation}\label{eq10}
\{x_1, x_2, \cdots , x_n, \theta , x_{n+1}, \cdots , x_m\}\,\equiv \,0\,(\,\mbox{mod}\, \mathcal{H}') \quad\mbox{for $0\le n\le m$.}
\end{equation}

It follows from Chon's Theorem that (\ref{eq9}) holds. Thus, the only thing we need to prove is (\ref{eq10}). We will prove (\ref{eq10}) by induction on $m$.

\medskip
Clearly, we have
\begin{equation}\label{eq11}
\{\theta \}=\theta, \quad \{x_1, \theta \}=\{ \theta , x_1 \}=\theta\bullet x_1 ,
\end{equation}
\begin{equation}\label{eq12}
\{\theta ,x_1, x_2\}=\{x_2, x_1, \theta \}=(\theta\bullet x_1 )\bullet x_2-
 (\theta\bullet x_2 )\bullet x_1 +\theta\bullet (x_1 \bullet x_2)
\end{equation}
and
\begin{equation}\label{eq13}
\{x_1, \theta , x_2\}=(\theta\bullet x_1 )\bullet x_2+
 (\theta\bullet x_2 )\bullet x_1 -\theta\bullet (x_1 \bullet x_2).
\end{equation}
It follows from (\ref{eq11}), (\ref{eq12}) and (\ref{eq13}) that (\ref{eq10}) holds for $m=0$, $1$ and $2$.

\medskip
We now assume that $m\ge 3$ and 
\begin{equation}\label{eq14}
\{z_1, \cdots , z_s, \theta , z_{s+1}, \cdots , z_t\}\,\equiv \,0\,(\,\mbox{mod}\, \mathcal{H}'), \end{equation}
where $0\le s\le t<m$, $\theta\in \Theta$ and $z_1$, $\cdots$, $z_t\in X$. If $0\le h< m-n$ and $m-n\ge 1$, then (\ref{eq9}) and (\ref{eq14}) imply
\begin{eqnarray*}
\mathcal{H}'&\ni& 
8\{x_1, \cdots , x_n, \theta , x_{n+1}, \cdots , x_{n+h}\}\bullet 
\{ x_{n+h+1}, \cdots , x_{m}\}\\
&=&2(x_1 \cdots  x_n \theta  x_{n+1} \cdots  x_{n+h}
+x_{n+h} \cdots  x_{n+1}\theta x_n \cdots  x_1)\bullet\\
&&\quad \bullet (x_{n+h+1} \cdots  x_{m}+x_{m} \cdots  x_{n+h+1})\\
&=&(x_1 \cdots  x_n \theta  x_{n+1} \cdots  x_{n+h}
+x_{n+h} \cdots  x_{n+1}\theta x_n \cdots  x_1)\cdot\\
&&\quad \cdot (x_{n+h+1} \cdots  x_{m}+x_{m} \cdots  x_{n+h+1})+\\
&&\quad +(x_{n+h+1} \cdots  x_{m}+x_{m} \cdots  x_{n+h+1})\cdot\\
&&\quad \cdot (x_1 \cdots  x_n \theta  x_{n+1} \cdots  x_{n+h}
+x_{n+h} \cdots  x_{n+1}\theta x_n \cdots  x_1)\\
&=&\underbrace{x_1 \cdots  x_n \theta  x_{n+1} \cdots  x_{m}}_1+
\underbrace{x_1 \cdots  x_n \theta  x_{n+1} \cdots  x_{n+h}x_{m} \cdots  x_{n+h+1}}_2+\\
&&\quad +\underbrace{x_{n+h} \cdots  x_{n+1}\theta x_n \cdots  x_1x_{n+h+1} \cdots  x_{m}}_3+\\
&&\quad +\underbrace{x_{n+h} \cdots  x_{n+1}\theta x_n \cdots  x_1x_{m} \cdots  x_{n+h+1}}_4+\\
&&\quad +\underbrace{x_{n+h+1} \cdots  x_{m}x_1 \cdots  x_n \theta  x_{n+1} \cdots  x_{n+h}}_4+\\
&&\quad +\underbrace{x_{n+h+1} \cdots  x_{m}x_{n+h} \cdots  x_{n+1}\theta x_n \cdots  x_1}_2+\\
&&\quad +\underbrace{x_{m} \cdots  x_{n+h+1}x_1 \cdots  x_n \theta  x_{n+1} \cdots  x_{n+h}}_3
+\underbrace{x_{m} \cdots  x_{n+1}\theta x_n \cdots  x_1}_1\\
&=&2\{x_1, \cdots , x_n, \theta ,x_{n+1} ,\cdots , x_{m}\}+\\
&&\quad +2\{x_1, \cdots , x_n, \theta , x_{n+1}, \cdots , x_{n+h}, x_{m}, \cdots , x_{n+h+1}\}+\\
&&\quad +2\{x_{n+h}, \cdots , x_{n+1}, \theta , x_n, \cdots , x_1, x_{n+h+1}, \cdots , x_{m}\}+\\
&&\quad +2\{x_{n+h}, \cdots , x_{n+1}, \theta , x_n, \cdots , x_1, x_{m}, \cdots , x_{n+h+1}\}
\end{eqnarray*}
or
\begin{eqnarray}\label{eq15}
&&\{x_1, \cdots , x_n, \theta ,x_{n+1} ,\cdots , x_{m}\}+\nonumber\\
&&\quad +\{x_1, \cdots , x_n, \theta , x_{n+1}, \cdots , x_{n+h}, x_{m}, \cdots , x_{n+h+1}\}+\nonumber\\
&&\quad +\{x_{n+h}, \cdots , x_{n+1}, \theta , x_n, \cdots , x_1, x_{n+h+1}, \cdots , x_{m}\}+\nonumber\\
&&\quad +\{x_{n+h}, \cdots , x_{n+1}, \theta , x_n, \cdots , x_1, x_{m}, \cdots , x_{n+h+1}\}\nonumber\\
&\,\equiv \,&0\,(\,\mbox{mod}\, \mathcal{H}') \quad\mbox{for $0\le h< m-n$ and $m-n\ge 1$.}
\end{eqnarray}

Let $h=n=0$ in (\ref{eq15}), we get
\begin{eqnarray*}
&&\{\theta , x_1, \cdots , x_m\}+\{\theta , x_m, \cdots , x_1\}+\\
&&\quad +\{\theta , x_1, \cdots , x_m\}+\{\theta , x_m, \cdots , x_1\}\,\equiv \,0\,(\,\mbox{mod}\, \mathcal{H}')
\end{eqnarray*}
or
\begin{equation}\label{eq16}
\{\theta , x_1, x_2, \cdots , x_m\}
\,\equiv \, -\{\theta , x_m, \cdots , x_2, x_1\}\,(\,\mbox{mod}\, \mathcal{H}'). 
\end{equation}

Let $h=m-n-1$ in (\ref{eq15}), we get
\begin{eqnarray*}
&&\{x_1, \cdots , x_n, \theta ,x_{n+1} ,\cdots , x_{m}\} +\{x_1, \cdots , x_n, \theta , x_{n+1}, \cdots , x_{m-1}, x_{m}\}+\nonumber\\
&&\quad +\{x_{m-1}, \cdots , x_{n+1}, \theta , x_n, \cdots , x_1, x_{m}\} +\nonumber\\
&&\quad +\{x_{m-1}, \cdots , x_{n+1}, \theta , x_n, \cdots , x_1, x_{m}\}
\,\equiv \, 0\,(\,\mbox{mod}\, \mathcal{H}')
\end{eqnarray*}
or
\begin{eqnarray}\label{eq17}
&&\{x_1, \cdots , x_n, \theta ,x_{n+1} ,\cdots , x_{m}\}\,\equiv \, 
-\{x_{m-1}, \cdots , x_{n+1}, \theta , x_n, \cdots , x_1, x_{m}\} \nonumber\\
&\,\equiv \,& -\{x_{m}, x_1, \cdots , x_{n}, \theta , x_{n+1}, \cdots ,  x_{m-1}\} \,(\,\mbox{mod}\, \mathcal{H}').
\end{eqnarray}

It follows that
\begin{eqnarray}
&&\{x_1, \cdots , x_n, \theta ,x_{n+1} ,\cdots , x_{m}\}\stackrel{(\ref{eq17})}{\,\equiv \,} 
 -\{x_{m}, x_1, \cdots , x_{n}, \theta , x_{n+1}, \cdots ,  x_{m-1}\}\nonumber \\
&\stackrel{(\ref{eq17})}{\,\equiv \,}&
(-1)^2\{x_{m-1}, x_m, x_1, \cdots , x_{n}, \theta , x_{n+1}, \cdots ,  x_{m-2}\}
\stackrel{(\ref{eq17})}{\,\equiv \,}\cdots \stackrel{(\ref{eq17})}{\,\equiv \,}\nonumber \\
\label{eq18}&\stackrel{(\ref{eq17})}{\,\equiv \,}&
(-1)^{m-n}\{x_{n+1}, \cdots , x_m, x_1, \cdots , x_{n}, \theta \}
\,(\,\mbox{mod}\, \mathcal{H}')\\
\label{eq19}&\,\equiv \,&(-1)^{m-n}\{\theta , x_n, \cdots , x_1, x_m, \cdots , x_{n+1}\}
\,(\,\mbox{mod}\, \mathcal{H}')\\
&\stackrel{(\ref{eq16})}{\,\equiv \,}&(-1)^{m-n+1}\{\theta , x_{n+1}, \cdots , x_m, x_1\cdots , x_{n}\} \nonumber \\
&\stackrel{(\ref{eq17})}{\,\equiv \,}&(-1)^{m-n+2}\{x_n, \theta , x_{n+1}, \cdots , x_m, x_1\cdots , x_{n-1}\}\stackrel{(\ref{eq17})}{\,\equiv \,}\cdots \stackrel{(\ref{eq17})}{\,\equiv \,} \nonumber \\
&\stackrel{(\ref{eq17})}{\,\equiv \,}&(-1)^{m-n+(n+1)}\{x_1, \cdots , x_n, \theta , x_{n+1}, \cdots , x_m\} \nonumber \\
\label{eq20}
&=&(-1)^{m+1}\{x_1, \cdots , x_n, \theta , x_{n+1}, \cdots , x_m\}\,(\,\mbox{mod}\, \mathcal{H}').
\end{eqnarray}

By (\ref{eq20}), (\ref{eq10}) holds for even $m$. From now on assume that $m$ is odd and 
$m\ge 3$. By (\ref{eq19}), to complete the proof of Proposition~\ref{pr3.1}, the only remaining thing we need to prove is
\begin{equation}\label{eq21}
\{\theta , x_1, \cdots , x_m\}\,\equiv \,0\,(\,\mbox{mod}\, \mathcal{H}') 
\quad\mbox{for $m$ is odd and $m\ge 3$.}
\end{equation}

Letting $n=0$ and $h=m-2$ in (\ref{eq15}), we get
\begin{eqnarray}\label{eq22}
&&\{ \theta ,x_{1} ,\cdots , x_{m}\} 
+\{\theta , x_{1}, \cdots , x_{m-2}, x_{m},  x_{m-1}\}+\nonumber\\
&&\quad +\{x_{m-2}, \cdots , x_{1}, \theta , x_{m-1}, x_{m}\}+\nonumber\\
&&\quad +\{x_{m-2}, \cdots , x_{1}, \theta , x_{m},  x_{m-1}\}
\,\equiv \,0\,(\,\mbox{mod}\, \mathcal{H}') 
\end{eqnarray}

It follows from (\ref{eq18}) and (\ref{eq22}) that
\begin{equation}\label{eq23}
\{ \theta ,x_{1} ,\cdots , x_{m}\} \,\equiv \,
-\{\theta , x_{1}, \cdots , x_{m-2}, x_{m},  x_{m-1}\}(\,\mbox{mod}\, \mathcal{H}'). 
\end{equation}

By (\ref{eq16}), we have
\begin{equation}\label{eq24}
\{ \theta ,x_{1} ,\cdots , x_{m}\} \,\equiv \,
-\{ x_{1}, \cdots , x_{m-1}, x_{m}, \theta \}(\,\mbox{mod}\, \mathcal{H}'). 
\end{equation}

Since the cycle $\left(\begin{array}{cccccc}\theta&x_1&x_2&\cdots&x_{m-1}&x_m\\
x_1&x_2&x_3&\cdots&x_m&\theta\end{array}\right)$ and the transposition
$\left(\begin{array}{cccccc}\theta&x_1&\cdots&x_{m-2}&x_{m-1}&x_m\\
\theta&x_1&\cdots&x_{m-2}&x_{m}&x_{m-1}\end{array}\right)$ generate the symmetric group $\mathbf{S}_{m+1}$ on the set $\{\theta, x_1, \cdots , x_m\}$, it follows from (\ref{eq23}) and (\ref{eq24}) that
\begin{equation}\label{eq25}
\{ \theta ,x_{1} ,\cdots , x_{m}\} \,\equiv \,
(\mbox{sign}\,\sigma)\{ \sigma (\theta), \sigma (x_{1}), \cdots , \sigma (x_{m}) \}(\,\mbox{mod}\, \mathcal{H}')
\end{equation}
for any permutation $\sigma =\left(\begin{array}{ccccc}\theta&x_1&x_2&\cdots&x_m\\
\sigma (\theta)&\sigma (x_1)&\sigma (x_2)&\cdots&\sigma (x_m)\end{array}\right)$
in $\mathbf{S}_{m+1}$, where $\mbox{sign}\,\sigma$ denotes the sign of the permutation $\sigma$.

\medskip
If $m=3$, then there is a permutation $\sigma =\left(\begin{array}{cccc}\theta&x_1&x_2&x_3\\
\theta&\sigma (x_1)&\sigma (x_2)&\sigma (x_3)\end{array}\right)$ on the set
$\{\theta, x_1, x_2, x_3\}$ such that $\sigma (x_1)<\sigma (x_2)<\sigma (x_3)$. Using this permutation $\sigma$ and (\ref{eq25}), we get
$$
\{ \theta ,x_{1}, x_2 , x_{3}\} \,\equiv \,
\pm \{ \theta, \sigma (x_{1}), \sigma (x_{2}) , \sigma (x_{3}) \}
\,\equiv \, 0(\,\mbox{mod}\, \mathcal{H}'),
$$
which proves that (\ref{eq21}) holds for $m=3$.

\medskip
In the remaining of the proof, we assume that $m$ is odd and $m\ge 5$. Let $n=0$ and $h=3$ in 
(\ref{eq15}), we get
\begin{eqnarray}\label{eq26}
&&\{ \theta ,x_{1} ,\cdots , x_{m}\}
+\{ \theta , x_{1}, x_2 , x_{3}, x_{m}, \cdots , x_5, x_{4}\} +\nonumber\\
&&\quad +\{x_{3}, x_2 , x_{1}, \theta , x_{4}, x_5, \cdots , x_{m}\}+\nonumber\\
&&\quad
+\{x_{3}, x_2 , x_{1}, \theta ,  x_{m}, \cdots , x_5, x_{4}\}\,\equiv \,  0\,(\,\mbox{mod}\, \mathcal{H}').
\end{eqnarray}

Let $\tau$ and $\mu$ be two permutations on the set $\{\theta, x_1, x_2, \cdots , x_m\}$ defined by
$$
\tau :=\left(\begin{array}{cccccc}\theta&x_1&x_2&\cdots&x_{m-1}&x_m\\
x_1&x_2&x_{3}&\cdots&x_m&\theta\end{array}\right)
$$
and
$$
\mu :=\left(\begin{array}{cccccccc}x_3&x_2&x_1&\theta&x_m&x_{m-1}&\cdots&x_4\\
\theta&x_1&x_2&x_3&x_m&x_{m-1}&\cdots&x_4\end{array}\right).
$$
Then $\mu$ is even and 
$$
\tau ^4 =\left(\begin{array}{cccccccc}\theta&x_1&\cdots&x_{m-4}&x_{m-3}&x_{m-2}&x_{m-1}&x_m\\
x_4&x_5&\cdots&x_m&\theta&x_1&x_2&x_3\end{array}\right).
$$

Using $\tau$, $\mu$ and (\ref{eq25}), (\ref{eq26}) implies
$$
4\{ \theta ,x_{1} ,\cdots , x_{m}\}\,\equiv \,  0\,(\,\mbox{mod}\, \mathcal{H}'),
$$
which proves that (\ref{eq21}) holds for odd $m\ge 5$.

\medskip
This completes the proof of Proposition~\ref{pr3.1}.

\hfill\raisebox{1mm}{\framebox[2mm]{}}


\bigskip
\begin{thebibliography}{9}
\bibitem{N} Nathan Jacobson, \textsl{Structure and Representations of Jordan Algebras}, Amer. Math. Soc. Colloq. Pub., Vol.39, Amer. Math. Soc., Providence, R.I., 1968.
\bibitem{Liu} Keqin Liu \textsl{Generalizations of Jordan algebras and Malcev algebras},
arXiv:math.RA/0606628 v1, 24 Jun 2006.
\bibitem{Liu3} Keqin Liu, \textsl{Hu-Liu Products}, Research Monographs in Mathematics {\bf 3}, 153 Publishing ({\it to appear}).
\bibitem{M} Kevin McCrimmon, \textsl{A Taste of Jordan Algebras}, Univertext, Springer, 2004.
\end{thebibliography}
\end{document}